\documentclass{article}

\author{Claire Levaillant\footnote{This research grew during a visit of the author at the Centro di Ricerca Matematica
Ennio de Giorgi of the Scuola Normale Superiore di Pisa. The author
thanks the center for its hospitality.}}
\title{Tangles of type $E_n$ and a reducibility criterion for the Cohen-Wales representation of the Artin
group of type $E_6$}

\usepackage{amsmath, amssymb}
\usepackage{amsthm}
\usepackage{palatino}
\usepackage{amscd}
\usepackage{graphicx}
\usepackage{epsfig}
\usepackage{verbatim}

\newcommand{\He}{\mathcal{H}}

\newcommand{\n}{\nu}
\newcommand{\Q}{\mathbb{Q}}

\newcommand{\al}{\alpha}

\newcommand{\ovl}{\overline}
\newcommand{\unsurr}{\frac{1}{r}}
\newcommand{\unsur}{\frac{1}}

\newcommand{\U}{\mathcal{U}}

\newcommand{\lb}{\lbrace}
\newcommand{\rb}{\rbrace}
\newcommand{\noin}{\noindent}

\newcommand{\LK}{\mathcal{L}\mathcal{K}}
\newcommand{\A}{\mathcal{A}}

\newcommand{\wh}{\widehat}

\newcommand{\ovs}{\overset}
\newcommand{\tpp}{\triangle^{++}}
\newcommand{\tp}{\triangle^{+}}
\newcommand{\tppp}{\triangle^{+++}}

\newcommand{\lr}{\ovs{\triangle^{l}}{w_{56}}}
\newcommand{\ta}{\ovs{\tppp}{w_{56}}}

\newcommand{\bl}{\blacklozenge}
\newcommand{\ab}{\mathcal{K}_6}

\begin{document}
\maketitle


\noin\textbf{Abstract.} We introduce tangles of type $E_n$ and
construct a representation of the Birman-Murakami-Wenzl algebra (BMW
algebra) of type $E_6$. As a representation of the Artin group of
type $E_6$, this representation is equivalent to the faithful
Cohen-Wales representation of type $E_6$ that was used to show the
linearity of the Artin group of type $E_6$. We find a reducibility
criterion for this representation and complex values of the
parameters for which the algebra is not semisimple. 

\section{Introduction and Main Results}
Birman-Murakami-Wenzl algebras were introduced in \cite{BIR} and
independently in \cite{MUR} in order to study the linearity of the
braid groups. The BMW algebra has the same generators as those of
the braid group and is a deformation of the Brauer centralizer
algebra. Like the Brauer algebra, the BMW algebra has a diagrammatic
version. It was shown by Morton and Wassermann in \cite{MOR} to be
isomorphic to the Morton-Traczyk algebra, an algebra of tangles.
Moreover, the BMW algebra is equipped with a trace functional such
that the Kauffman polynomial invariant of links is, after an
appropriate renormalization that trace (see \cite{BIR}). In
\cite{CGW}, Cohen, Gijsbers and Wales generalize the BMW algebra to
the simply laced Coxeter types $D$ and $E$ and in \cite{BAL} they
show that the BMW algebra of type $D_n$ is isomorphic to a tangle
algebra which they define. The generalized BMW algebra has the same
generators as those of the Artin group in the same way the original
BMW algebra has the same generators as those of the braid group. In
\cite{CL}, we use the tangle algebra defined in \cite{CGW} to
construct a representation of the BMW algebra of type $D_n$
$BMW(D_n)$, which as a representation of the Artin group $\A(D_n)$
is equivalent to the generalized Lawrence--Krammer representation
$\LK(D_n)$ introduced by Arjeh Cohen and David Wales in \cite{CW}.
We use this representation to deduce a reducibility criterion for
$\LK(D_n)$ as well as a conjecture which gives a criterion of
semisimplicity for $BMW(D_n)$. In this paper we propose tangles of
type $E_6$, and more generally $E_n$. We use them as a tool to
construct a representation of the BMW algebra of type $E_6$ (denoted
$BMW(E_6)$). As a representation of the Artin group of type $E_6$,
this representation is equivalent to the generalized
Lawrence--Krammer representation of type $E_6$ (denoted $\LK(E_6)$).
We use our representation to give a reducibility criterion for
$\LK(E_6)$ and derive some complex values of the parameters for
which
$BMW(E_6)$ is not semisimple.\\
In the following theorem, the representation $\LK(E_6)$ renamed
after its authors the Cohen-Wales representation of the Artin group
of type $E_6$ with parameters $t$ and $r$ is the one defined in
\cite{CW}.
\newtheorem{Theorem}{Theorem}
\begin{Theorem}
The Cohen-Wales representation of the Artin group of type $E_6$ with
non-zero complex parameters $t$ and $r$ is irreducible, except when
$$t\in\lb 1,-1,r^6,-r^{12},r^{24}\rb,$$
when it is reducible. Moreover, \\
\begin{center}
when $t=1$, there exists an invariant subspace of dimension $15$. \\
When $t=-1$, there exists an invariant subspace of dimension $30$.\\
When $t=r^6$, there exists an invariant subspace of dimension
$20$.\\
When $t=-r^{12}$, there exists an invariant subspace of dimension
$6$.\\
When $t=r^{24}$, there exists an invariant subspace of dimension
$1$.
\end{center}
\end{Theorem}
\begin{Theorem}
Let $l$, $m$ and $r$ be three non-zero complex numbers with
$m=\unsurr-r$. \\
When $$l\in\bigg\lb
r^3,-r^3,\unsur{r^3},-\unsur{r^3},-\unsur{r^9},r^9,\unsur{r^{21}},-r^{21}\bigg\rb$$
the BMW algebra of type $E_6$ with non-zero complex parameters $l$
and $m$ is not semisimple over the field $\Q(l,r)$.
\end{Theorem}
\indent The interest of the paper is double. The Lawrence-Krammer
representations have become of interest since the original
Lawrence--Krammer representation was used by Stephen Bigelow
\cite{BIG} and independently Daan Krammer \cite{KR} to show the long
open problem of the linearity of the braid group. The generalized
Lawrence--Krammer representations of types $D$ and $E$ introduced by
Cohen and Wales were in turn used by them to show the linearity of
the Artin groups of these types (see \cite{CW}). Other constructions
and proof of linearity can also be found in \cite{DIG}. Reducibility
criteria for these representations exist in types $A$ and $D$ (see
\cite{CRAS} and \cite{CL} respectively). In each case, reducibility
is shown for some complex specializations of the two parameters of
the representation, while the representations are generically
irreducible (see \cite{CGW}, \cite{MAR}, \cite{ZINN}, \cite{CRAS} in
type $A$ and \cite{CGW} and \cite{CL} in type $D$). More work allows
to determine the complete structure of the representation. This was
done in type $A$ in \cite{RS} and independently in \cite{PROC} by
different means; it was later achieved in type $D$ in \cite{CLA}.
The specializations for which the representation becomes reducible
give non-semisimplicity for the BMW algebra whose parameters are
related to the parameters of the representation. This is another
important aspect. There has been a long interest in studying the
semisimplicity of the BMW algebras. The problem in type $A$
originated in Result $(a)$ of \cite{WEN}, was further developed in
Theorem $2$ of \cite{CLDBW} and was completely solved in Theorem $B$
of \cite{RUI}, building upon \cite{ENY}. The problem in type $D$ is
stated as a conjecture in \cite{CL} and is to this date not solved. \\

We end this introduction by recalling the defining relations of
$BMW(E_6)$. First, we shall draw a Dynkin diagram of type
$E_6$.\vspace{-0.1cm}
\begin{center}\epsfig{file=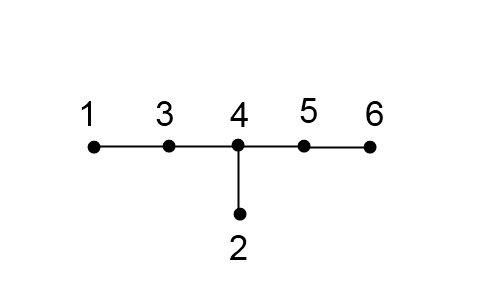, height=5cm}\end{center}
We write $i\sim j$ if nodes $i$ and $j$ are adjacent on the Dynkin
diagram and $i\not\sim j$ if nodes $i$ and $j$ are not adjacent on
the Dynkin diagram. The BMW algebra of type $E_6$ with nonzero
complex parameters $l$ and $m$ is an algebra over $\Q(l,m)$ with six
generators called the $g_i$'s, $1\leq i\leq 6$. It contains other
elements, namely the $e_i$'s, $1\leq i\leq 6$ that are related to
the $g_i$'s as below.
$$\begin{array}{ccccc}
(A1)&g_ig_j&=&g_jg_i&\text{if $i\not\sim j$}\\
(A2)&g_ig_jg_i&=&g_jg_ig_j&\text{if $i\sim j$}\\
(P)&e_i&=&\frac{l}{m}(g_i^2+m\,g_i-1)&\text{for all $i$}\\
(DL1)&g_ie_i&=&\unsur{l}\,e_i&\text{for all $i$}\\
(DL2)&e_ig_je_i&=&l\,e_i&\text{when $i\sim j$}
\end{array}$$

\noin $(A_1)$ and $(A_2)$ are the Artin group relations; $(P)$
defines each $e_i$ as a polynomial in $g_i$; $(DL1)$ and $(DL2)$ are
called the delooping relations because they are the algebraic
versions of the delooping relations on the tangles (see
$\S\,2$, point $(iii)$ of Definition $2.6$ of \cite{BAL}).\\
Some consequences of these defining relations are the following.
$$\begin{array}{ccccccc}
(I)& e_i^2&=&\delta\,e_i&&&\text{for all $i$ where
$\delta=1-\frac{l-\unsur{l}}{m}$}\\
(MA)& g_ig_je_i&=&e_je_i&=&e_jg_ig_j&\text{when $i\sim j$}\\
(R)&e_ie_je_i&=&e_i&&&\text{when $i\sim j$}
\end{array}$$

\noin $(I)$ expresses the fact that $\delta^{-1}\,e_i$ is an
idempotent; We call equalities in $(MA)$ "mixed Artin relations";
$(R)$ can be viewed as a
reduction. \\
A result of \cite{CWE} states that $BMW(E_6)$ is semisimple and free
over the integral domain $\mathbb{Z}[\delta,
\delta^{-1},l,l^{-1},m]/(m(1-\delta)-(l-\unsur{l}))$ of rank $1,\,
440,\, 585$ and that it is cellular in the sense of Graham and
Lehrer \cite{LEH} over suitable rings.
\section{Tangles of type $E_n$}

We use tangles as a tool to construct a representation of
$BMW(E_6)$. Here is how the algebra elements $e_i$'s and $g_i$'s are
represented in terms of tangles. \epsfig{file=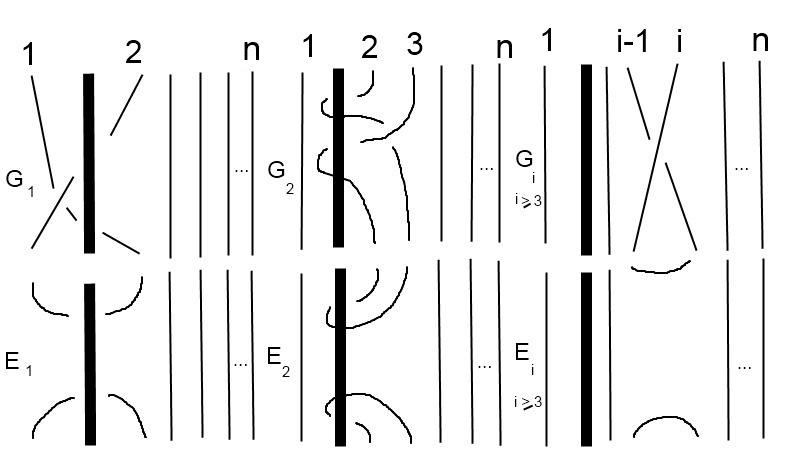,
height=8cm}\\
The rigid bar is called the pole. The strands can twist around the
pole. A strand twists around the pole if it goes behind the pole and
then over it, or the converse. An important aspect is that for a
strand twisting around the pole, the following relations are
satisfied.
\begin{center} \epsfig{file=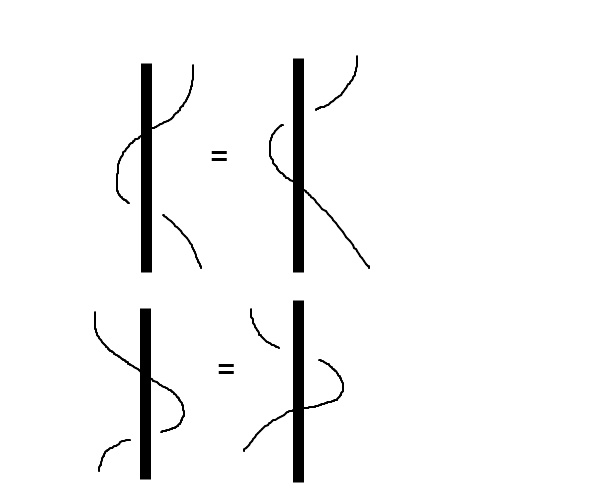, height=7cm}
\end{center}

\noin They are called the twist relations. We impose that the
tangles satisfy other relations, some of which are given in
\cite{BAL}. We add to the classical Reidemeister's moves
\newcounter{c}\setcounter{c}{2} \Roman{c} and \setcounter{c}{3}
\Roman{c} two pole-related Reidemeister's moves that we call
Reidemeister's moves of types \setcounter{c}{4} \Roman{c} and
\setcounter{c}{5} \Roman{c}. They will be respectively denoted by
\setcounter{c}{4}(R \Roman{c}) and \setcounter{c}{5} (R
\Roman{c}).\\
\epsfig{file=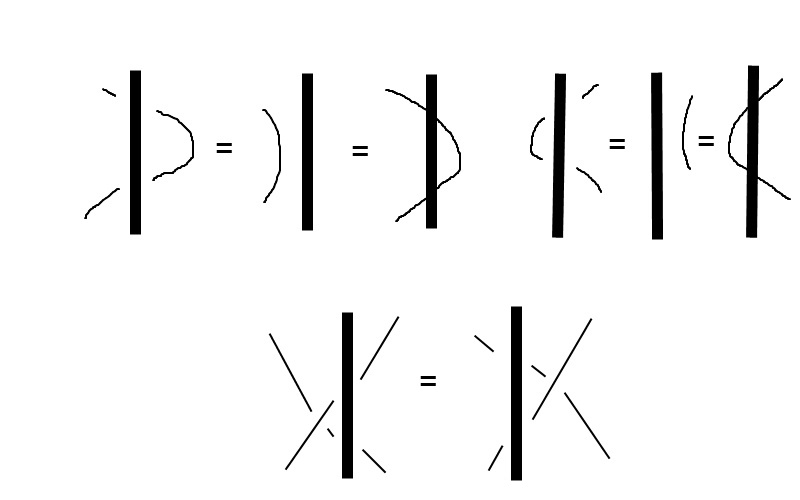, height=9cm}\\
\begin{center} \textit{The pole-related Reidemeister's
moves}\end{center}  \noin We now recall below the commuting relation
of \cite{BAL}.
\begin{center}
\epsfig{file=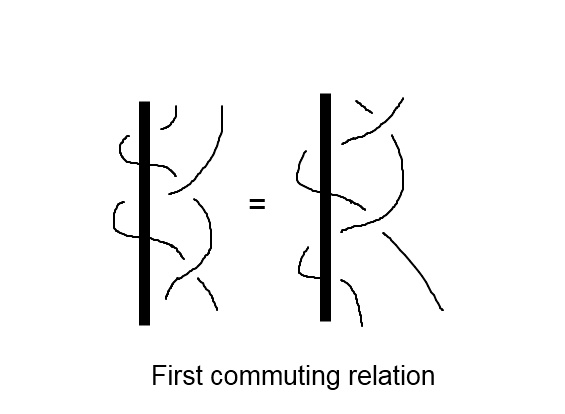, height=7cm}
\end{center}
We add a new commuting relation which is specific to type $E_6$ and
which is the following. \vspace{-0.45cm}\begin{center}
\epsfig{file=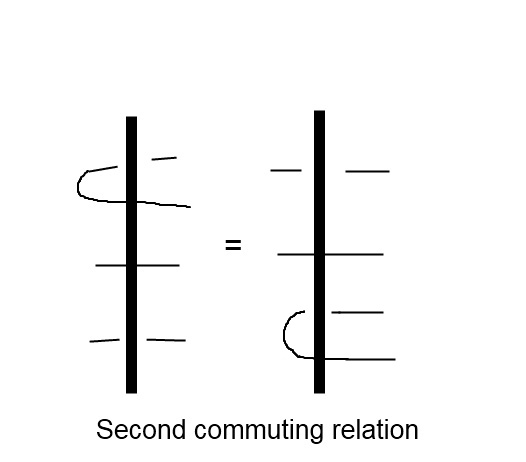, height=7cm}\end{center} We call it "commuting
relation of the second type" and we refer to it as (C
\setcounter{c}{2}\Roman{c}). We call the commuting relation of
\cite{BAL} "commuting relation of the first kind" and we will refer
to it as (C \setcounter{c}{1}\Roman{c}). It is a straightforward
verification to check that by using the twist relation, (R
\setcounter{c}{4}\Roman{c}) and (C \setcounter{c}{2}\Roman{c}), the
following relations are satisfied on the tangles
\begin{eqnarray*}
G_1\,E_2&=&E_2\,G_1\\
E_1\,G_2&=&G_2\,E_1\\
G_1\,G_2&=&G_2\,G_1\\
E_1\,E_2&=&E_2\,E_1
\end{eqnarray*}
We do the first one and leave the other ones to the reader.\\
\epsfig{file=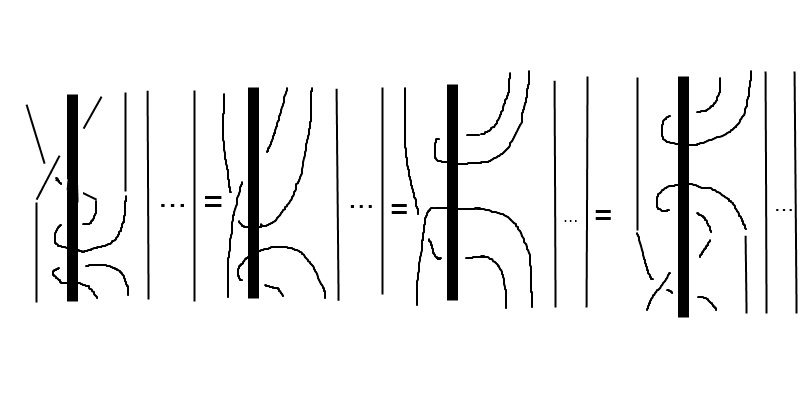, height=6cm}\\
\noin The first and third equalities are obtained by using (R
\setcounter{c}{4}\Roman{c}); the second equality is obtained by
using (C\setcounter{c}{2}\Roman{c}) in the following way. Divide the
middle tangle to the left into three vertical strips. The strip in
the middle contains the twist, the horizontal over-crossing and the
horizontal under-crossing. The strands in the right hand side strip
intersect the border of the middle strip in four points. The strands
in the left hand side strip intersect the border of the middle strip
in two points. Apply (C\setcounter{c}{2}\Roman{c}) in the middle
strip with the extremities of the strands reaching the points
mentioned above in the order in which they appear, while
leaving the other strips unchanged. Obtain the middle tangle to the right.\\\\
\noin Some closed pole loop relations are given in Proposition
$2.10$ of \cite{BAL}. Here, we show an additional proposition.
\newtheorem{Proposition}{Proposition}
\begin{Proposition} The tangles of type $E_6$ satisfy the following closed pole loop
relations that we call "the fourth closed pole loop relations".
\begin{center}\epsfig{file=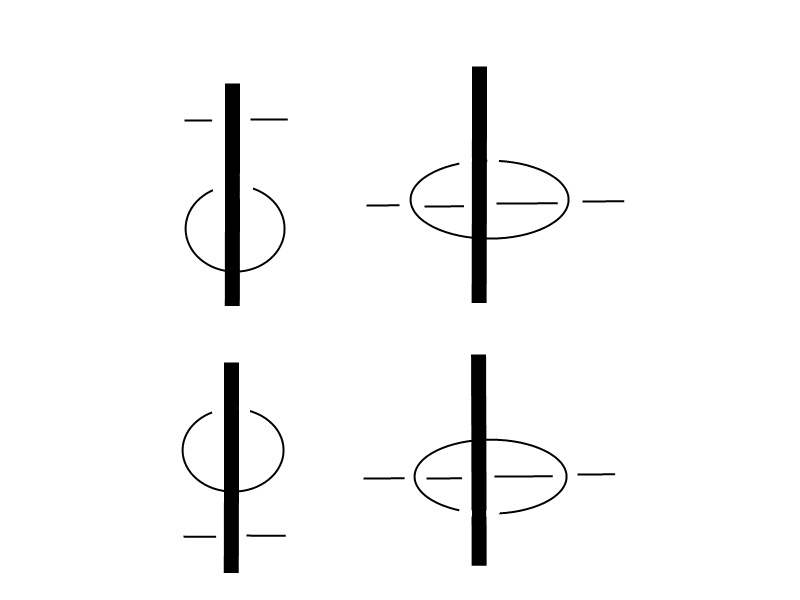,
height=7cm}\end{center} $\qquad\qquad\qquad\qquad$\textit{The fourth
closed pole loop relations}.
\end{Proposition}
\noin \textsc{Proof.} We partly use the same argument as in the
proof of Proposition $2.10$ of \cite{BAL}, point $(i)$. Indeed,
using (R
\newcounter{co}\setcounter{co}{2} \Roman{co}), we enlarge the closed
loop in such a way that it over-crosses the horizontal strand four
times. We then apply (R \setcounter{co}{5}\Roman{co}), followed by
(R \setcounter{co}{2}\Roman{co}). Contrary to \cite{BAL}, it is not
the commuting relation that plays a role here, but the new
Reidemeister's move (R \setcounter{co}{5}\Roman{co}). This property
is thus specific to tangles of type $E_6$. The moves are illustrated
on the figures below. \epsfig{file=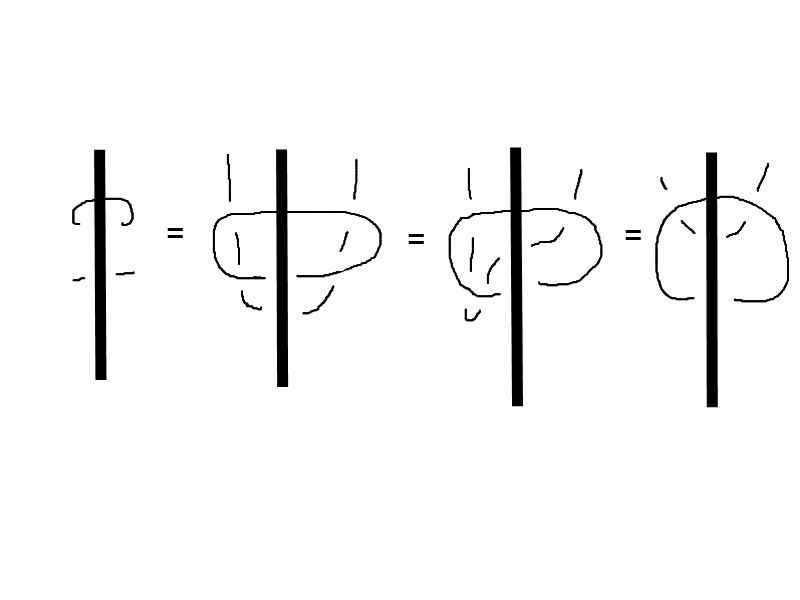,
height=9.5cm}\vspace{-4.2cm} \epsfig{file=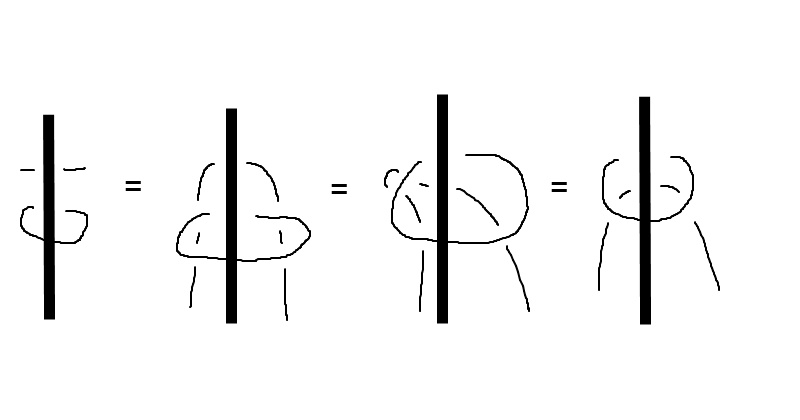, height=7.5cm}

\noin For all the other relations satisfied on the tangles, namely
the Kauffman skein relation, the self-intersection relations, the
idempotent relation, the first pole-related self-intersection
relation, the second pole-related self-intersection relation, the
first closed pole loop relation and their consequences, we refer the
reader to Definition $2.6$ and to Figures $6$, $9$, $13$ of
\cite{BAL}. With these relations, all the defining relations of the
BMW algebra are satisfied on the tangles $E_i$'s and $G_i$'s. Thus,
there is a morphism of algebras between the BMW algebra and the
tangle algebra, which sends $g_i$ to $G_i$ and $e_i$ to $E_i$. We do
not know whether this morphism is an isomorphism of algebras and do
not answer the question in this paper. This homomorphic image of the
BMW algebra in the tangle algebra is likely the whole tangle
algebra, but we do not consider this in our paper. We simply use the
tangles as a tool in order to get the actions in the representation
that we build. We then check with Maple that the map that we define
is indeed a representation of the BMW algebra.

\section{The representation}
We build a representation of $BMW(E_6)$ inside the vector space
$V_6$ over $\Q(l,r)$ spanned by vectors indexed by the $36$ positive
roots of a root system of type $E_6$. By definition, $r$ and
$-\unsurr$ are the two non-zero complex roots of the polynomial
$X^2+m\,X-1$, and so $m$ and $r$ are related by $m=\unsurr-r$. For
the vectors indexed by the positive roots, we introduce some
convenient notations. First, if we forget about node $1$ and the
edge joining nodes $1$ and $3$ on the Dynkin diagram, we obtain a
Dynkin diagram of type $D_5$ on nodes $2,3,4,5,6$. For the vectors
indexed by the positive roots issued from this diagram, we use the
same notations as in \cite{CL} with node $i+1$ now playing the role
of node $i$. So, for instance $\wh{w_{2j}}$, $j\geq 4$ (resp
$\wh{w_{23}}$) denotes the vector associated with the positive root
$\al_2+\al_4+\dots+\al_j$ (resp with the simple root $\al_2$); we
denote by $\wh{w_{3j}}$, $j\geq 4$ the vector associated with the
positive root $\al_2+\al_3+\al_4+\dots+\al_j$; we denote by
$\wh{w_{st}}$, $s\geq 4$ the vector associated with the positive
roots $\al_2+\al_3+2\,\al_4+\dots+2\,\al_s+\al_{s+1}+\dots+\al_t$.
Like we did in type $D$, we associate tangles to these vectors as on
the examples below. And so, to each of these vectors corresponds an
algebra element on which the $g_i$'s can act to the left.
\\\epsfig{file=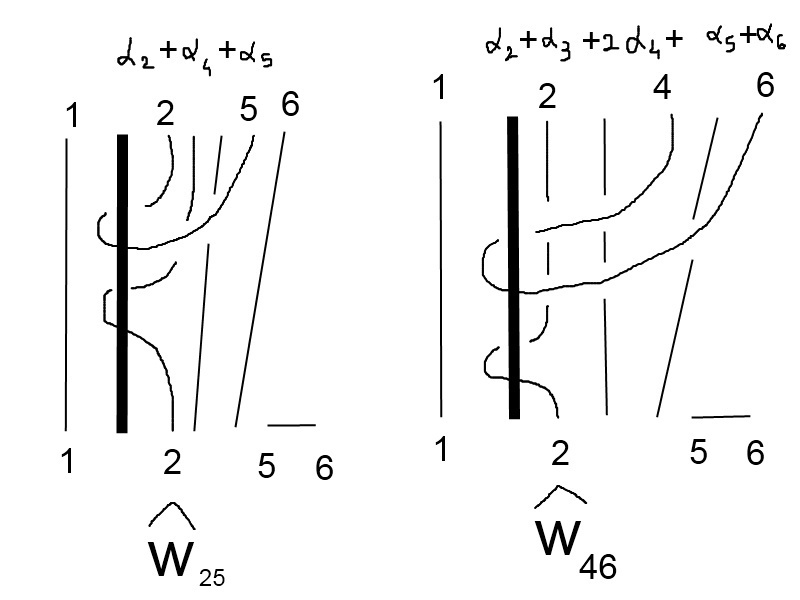, height=8cm}\\ In what
follows, by $\wh{w_{st}}$ we mean, depending on the context, the
algebra element or the tangle representing the algebra element or
the basis vector itself. By abuse of terminology, we will also
sometimes speak indifferently of the positive root, the vector
indexed by it, the algebra element associated with it or the tangle
representing the algebra element. Using these conventions, we now
present notations for the vectors indexed by a positive root whose
support contains node number $1$.
\newtheorem{Definition}{Definition}
\begin{Definition} We define
$$\begin{array}{cccc}
\ovs{\tp}{w_{st}}&=&g_1\,\wh{w_{st}}&\text{when $s\geq 3$}\\
\ovs{\tpp}{w_{st}}&=&g_3\,\ovs{\tp}{w_{st}}&\text{when $s\geq 4$}\\
\ta&=&g_4\,\ovs{\tpp}{w_{56}}\\
\lr&=&g_2\,\ta
\end{array}$$
\end{Definition}

\noin Following the notations of \cite{BOU}, we denote by

$$\begin{array}{ccccc}
a&c&d&e&f\\
&&b&&\end{array}$$

\noin the positive root
$\beta=a\,\al_1+b\al_2+c\al_3+d\al_4+e\al_5+f\al_6$.
And so, we see that\\\\

\noin A vector indexed by $\beta$ carries $\begin{cases} \text{a hat
if and only if $a=0$ and $b=1$}\\
\text{a triangle if and only
if $a=1$ and $b\neq 0$.}\\
\text{a triangle and an
additive sign if and only if $a=b=c=1$.}\\
\text{a triangle and two additive signs if and only if $c=d=2$} \\
\text{a triangle and three additive signs if and only if $b=1$ and
$d=3$}\\
\text{a triangle and an "l" as in "longest root" if and only if
$b=2$}
\end{cases}$\\\\

\noin Examples of tangles representing some elements of Definition
$1$ are given below, where we used the Bourbaki notation to refer to
the positive root on the top line.\\
\\\epsfig{file=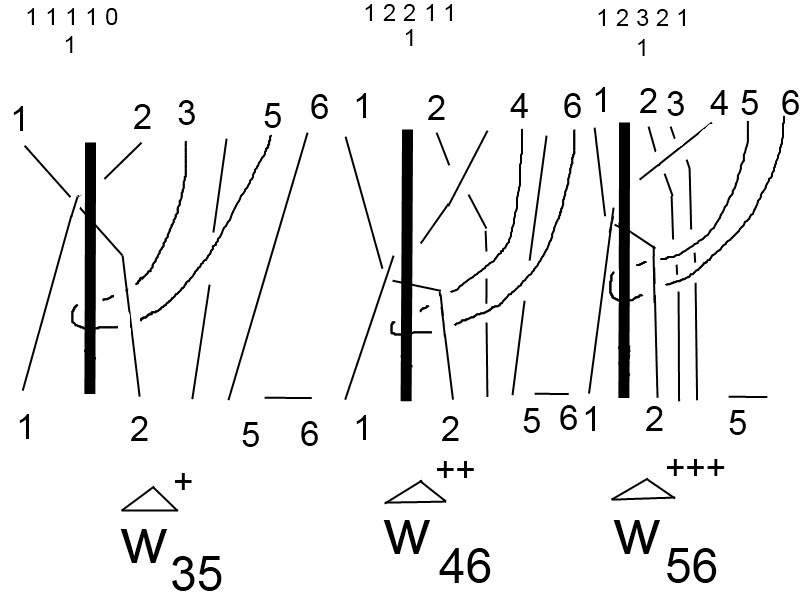, height=9cm}\\
\noin Note the number of additive signs gives the number of
crossings on the vertical strands. To obtain these diagrams from
those representing the $\wh{w_{ij}}$'s with $i\geq 3$, we use the
first commuting relation. This relation allows to switch the order
of the two strands that twist around the pole by putting one on top
of the other and conversely. The vertical strand that twists around
the pole used to under-cross the horizontal strand twisting around
the pole. After using the first commuting relation, it over-crosses
it. Changing the order of these
two strands then allows to use the pole-related Reidemeister's move of type \setcounter{co}{4} (\Roman{co}).  \\
For clarity, we now list in the table below the $36$ positive roots
of type $E_6$ and the vectors of $V_6$ indexed by them.

$$\begin{array}{cc}
\al_1&w_{12}\\
\al_1+\al_3+\dots+\al_t&w_{1t},\;3\leq t\leq 6\\
\al_{i+1}+\dots+\al_j&w_{ij},\;2\leq i< j\leq 6\\
\al_2&\wh{w_{23}}\\
\al_2+\al_4+\dots+\al_t&\wh{w_{2t}},\;4\leq t\leq 6\\
\al_2+\al_3+\al_4+\dots+\al_t&\wh{w_{3t}},\;t\geq 4\\
\al_2+\al_3+2\,\al_4+\dots+2\,\al_s+\al_{s+1}+\dots+\al_t&\wh{w_{st}},\;4\leq
s<t\leq 6\\
\al_1+\al_2+\al_3+\al_4+\dots+\al_t&\ovs{\tp}{w_{3t}},\;t\geq 4\\
\al_1+\al_2+\al_3+2\,\al_4+\dots+2\,\al_s+\al_{s+1}+\dots+\al_t&\ovs{\tp}{w_{st}},\;4\leq
s<t\leq 6\\
\al_1+\al_2+2\,\al_3+2\,\al_4+\dots+2\,\al_s+\al_{s+1}+\dots+\al_t&\ovs{\tpp}{w_{st}},\;4\leq
s<t\leq 6\\
\al_1+\al_2+2\,\al_3+3\,\al_4+2\,\al_5+\al_6&\ta\\
\al_1+2\,\al_2+2\,\al_3+3\,\al_4+2\,\al_5+\al_6&\lr
\end{array}$$

\noin The algebra elements associated with the positive roots are
gathered in the following table and in Definition $1$.
$$\begin{array}{cc}
e_i\dots e_6&w_{i-1,i},\;i\geq 3\\
e_1\,e_3\,e_4\,e_5\,e_6&w_{12}\\
g_j\dots g_{i+2}\,e_{i+1}\dots e_6&w_{ij},\;j\geq
i+2\;\text{and}\;i\neq 1\\
g_t\dots\,g_3\,e_1e_3e_4e_5e_6&w_{1t},\;t\geq 3\\
e_2e_4e_5e_6&\wh{w_{23}}\\
g_j\dots\,g_4\,e_2e_4e_5e_6&\wh{w_{2j}},\;j\geq 4\\
g_i\dots\,g_3g_j\dots\,g_4\,e_2e_4e_5e_6&\wh{w_{ij}},\;3\leq i<j\leq
6\end{array}$$

\noin Let $F$ denote the ground field $\Q(l,r)$ and $B(E_6)$ the BMW
algebra of type $E_6$ formerly denoted by $BMW(E_6)$. Our
representation is built inside the $B(E_6)$-module

$$B(E_6)e_6\bigg/\bigg(<B(E_6)e_ie_jB(E_6)>_{i\not\sim j}\cap\,
B(E_6)e_6\bigg)\otimes_{\He}F,$$ where $\He$ is the Hecke algebra of
type $A_5$ with generators $g_2$, $g_4$, $g_3$, $g_1$ and $z$. These
generators satisfy the relations $(A_1)$, $(A_2)$ from the defining
Artin relations of the BMW algebra, the relations
$g_i^2+m\,g_1=1\;\;\text{for all $i$}$, the relation $z^2+m\,z=1$
and the relations
$$\left\lb\begin{array}{cccc}z\,g_i&=&g_i\,z&\text{for all
$i\in\lb 2,3,4\rb$}\\
g_1\,z\,g_1&=&z\,g_1\,z&
\end{array}\right.$$

\noin We denote by $M$ the left $B(E_6)$-module
$$M=B(E_6)e_6\bigg/\bigg(<B(E_6)e_ie_jB(E_6)>_{i\not\sim j}\cap \,B(E_6)e_6\bigg)$$ We
provide $M$ with a structure of right $\He$-module by the following
actions.\\
The $g_i$'s act to the right of elements of $M$ by simply
multiplying them to the right by $g_i$ in $M$.\\
The Hecke algebra element $z$ acts to the right of an element of $M$
by multiplying it to the right in $M$ by
$$\xi=\unsur{\delta^2}\,e_6\,e_5\,e_4\,e_2\,g_3\,e_2\,e_4\,e_5\,e_6$$
The tangle representing $\xi$ is the following. \begin{center}
\epsfig{file=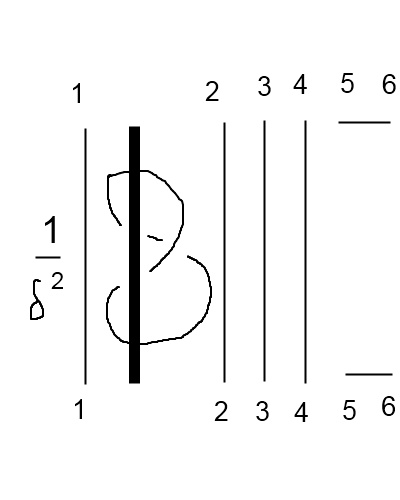, height=7.5cm}
\end{center}
Note that applying the commuting relation of the second type on
$\xi$ leaves the tangle invariant. It is a result of \cite{BAL} that
the $(0,0)$-tangle that has the shape of an eight, called $\Xi^{+}$
in \cite{BAL} commutes with any twist around the pole, as
illustrated on Figure $9$ of \cite{BAL}. Then, by making the tangle
$\Xi^{+}$ as small as we want, we see that $\Xi^{+}$ commutes with
$g_2$. Hence $\xi$ commutes with $g_2$. A quick glance at the
geometric definition of $\xi$ shows that $\xi$
also commutes with $g_3$ and $g_4$. 
We now show that $\xi$ and $g_1$ satisfy the following relation.
\begin{Proposition}
$g_1\,\xi\,g_1=\xi\,g_1\,\xi$
\end{Proposition}
\noin\textsc{Proof.} Follows from the relations $(I)$ and $(R)$ and from the Artin relation $g_1g_3g_1=g_3g_1g_3$. \\


\noin The rest of the proof that $M$ is a right $\He$-module for the
action provided above is up to an appropriate change in the indices
similar to the proof of Claim $1$ of \cite{CL}. Recall that
$r^2+m\,r=1$. We then provide the ground field $F$ with a structure
of left $\He$-module by letting the generators of $\He$ act on $1$
in the following way.
$$\left\lb\begin{array}{cccc} g_i.\,1&=&r&\text{for all $i$}\\
z.\,1&=&r&
\end{array}\right.$$
\noin We obtain a left representation of $BMW(E_6)$ by considering
the tensor product $M\otimes_{\He}F$. \\

We now use the tangles as a tool to compute the left actions of the
$g_i$'s on the basis vectors of $V_6$ inside the left
$BMW(E_6)$-module $M\otimes_{\He}F$. The fact that we work inside
$M\otimes_{\He}F$ allows to multiply a tangle at the bottom by $g_i$
(resp $g_i^{-1}$) for any $i\in\lb 1,2,3,4\rb$ at the cost of a
multiplication (resp division) by $r$. It also allows to replace a
pole related self intersection by a simple twist around the pole at
the cost of a factor $r$ or $\unsurr$, depending on the sign of the
crossing (see Figure $9$ of \cite{BAL} and the explanations in
$\S\,2.2$ of \cite{CL}). The move is illustrated on the figure
below. \\
\epsfig{file=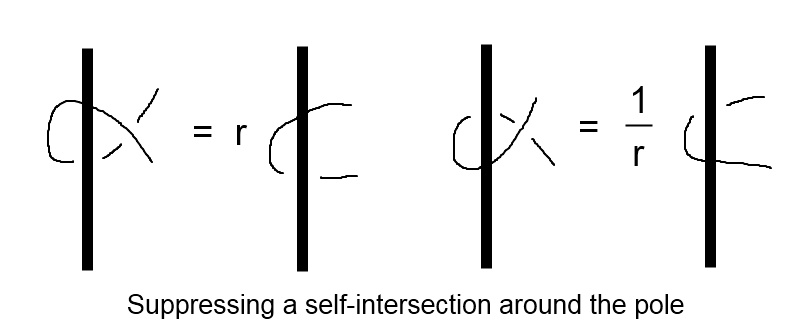, height=5cm}
\begin{center}
\textit{Fig $1$}\end{center} Finally it allows to suppress a
pole-related self intersection as on Figure $13$ of \cite{BAL} at
the cost of a multiplication by $r$. We succeed to define a
representation of $BMW(E_6)$ inside $V_6$. We define a map

$$\n^{(6)}:\begin{array}{ccc}BMW(E_6)&\longrightarrow &
\text{End}_{F}(V_6)\\
g_i&\longmapsto& \n_i
\end{array},$$

\noin with the endomorphisms $\n_i$'s given below.

\noin First, we give the left action by $g_1$ on our basis vectors.
When it not described, the action on the omitted basis vectors is a
multiplication by $r$.

\begin{eqnarray}
\n_1(w_{12})&=&\unsur{l}\,w_{12}\\
\forall j\geq 3,\,\qquad\;\n_1(w_{1j})&=&w_{2j}\\
\forall j\geq 3,\,\qquad\;\n_1(w_{2j})&=&w_{1j}+m\,r^{j-3}\,w_{12}-m\,w_{2j}\\
\forall j\geq 3,\,\qquad\;\n_1(\wh{w_{2j}})&=&r\,\wh{w_{2j}}\\
\forall s\geq 3,\,\qquad\;\n_1(\wh{w_{sj}})&=&\overset{\triangle^{+}}{w_{sj}}\\
\forall t\geq 4,\,\qquad\;\n_1(\ovs{\tp}{w_{3t}})&=&\wh{w_{3t}}-m\,\ovs{\tp}{w_{3t}}+\frac{m}{l}\,r^{t-6}\,w_{12}\\
\forall s\geq 4,\,\qquad\;\n_1(\ovs{\tp}{w_{sj}})&=&\wh{w_{sj}}-m\,\ovs{\tp}{w_{sj}}+\frac{m}{l}\,r^{s+t-9}\,w_{12}\\
\forall s\geq
4,\,\qquad\;\n_1(\;\ovs{\tpp}{w_{sj}})&=&r\,\;\ovs{\tpp}{w_{sj}}\\
\n_1(\;\ovs{\tppp}{w_{56}})&=&r\;\;\ovs{\tppp}{w_{56}}\\
\n_1(\;\lr)&=&r\;\lr
\end{eqnarray}

\noin Equality $(4)$ is obtained by using the second commuting
relation as on the figure below.\\
\epsfig{file=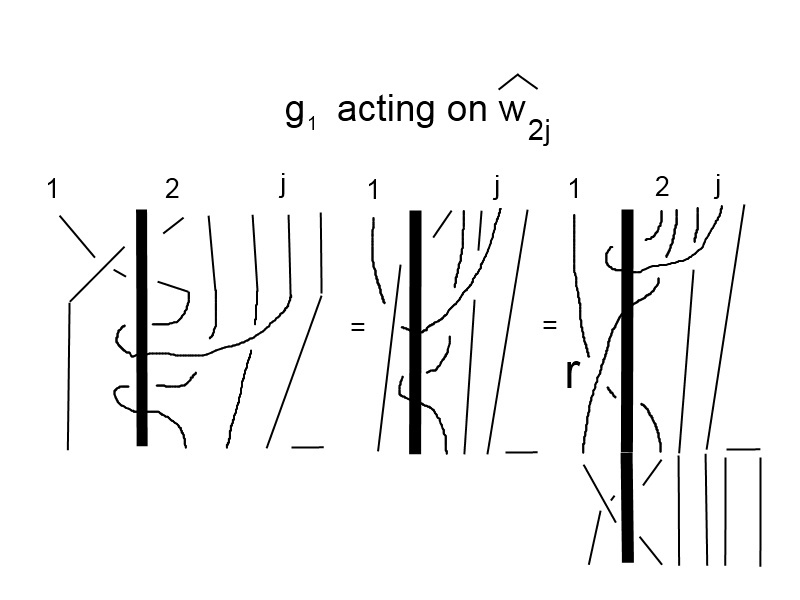, height=8.5cm}\\
\noin To derive equality $(6)$, we must compute the action of $e_1$
on $\wh{w_{3t}}$. To do so, we use the tangles. We apply the
commuting relation of the first kind and a Reidemeister's move of
type \newcounter{compteur}\setcounter{compteur}{4} \Roman{compteur}.
Then, we apply the commuting relation of the second kind. We
multiply the result at the bottom by $g_1$ at the cost of a division
by $r$. Then we apply (R\setcounter{compteur}{5} \Roman{compteur})
to move the crossing issued from $g_1$ to the right hand side of the
pole. Next, we successively apply (R\setcounter{compteur}{4}
\Roman{compteur}) and (R \setcounter{compteur}{2} \Roman{compteur}).
We then remove the crossings between the vertical strands at the
cost of a factor $r^{t-4}$. Finally, we apply (C
\setcounter{compteur}{2}\Roman{compteur}) a second time. When doing
so, we create a pole-related self-intersection using the terminology
of \cite{BAL}. As explained in \cite{CL} we can replace it with a
simple twist around the pole by dividing by a factor $r$. Hence the
term in $\frac{m}{l}r^{t-6}\,w_{12}$ in $(6)$. Some of the moves
have been gathered on the figure below.
\epsfig{file=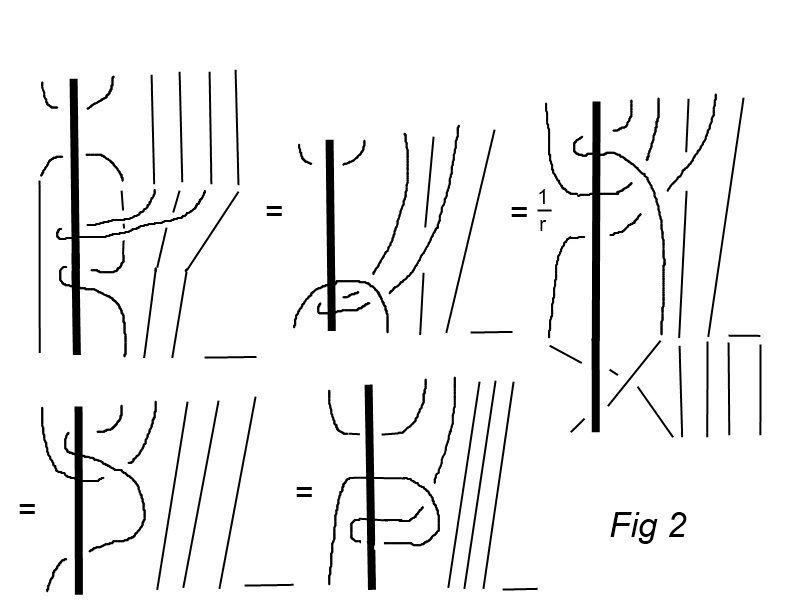, height=9cm}\\\\
\noin Equality $(7)$ is now obtained from Equality $(6)$ as follows.
\begin{eqnarray*}
\forall\,s\geq 4,\;\; e_1\wh{w_{st}}&=&e_1\,g_{s,4}\wh{w_{3t}}\\
&=&g_{s,4}\,e_1\wh{w_{3t}}\\
&=&r^{t-6}\,g_{s,4}\,w_{12}\\
&=&r^{s+t-9}\,w_{12}
\end{eqnarray*}

\noin Equalities $(8)$, $(9)$ and $(10)$ follow from the Artin group
relations and from the
fact that when $s\geq 4$, we have $g_3\wh{w_{st}}=r\,\wh{w_{st}}$.\\

The action by $g_2$ on the $\wh{w_{st}}$'s with $2\leq s\leq n-1$
and $3\leq t\leq n$ is given by reading the action
$\n_1(\wh{w_{s-1,t-1}})$ on the formulas $(1)$, $(2)$, $(3)$, $(6)$
and $(7)$ of Theorem $1$ of \cite{CL} and by incrementing all the
indices of the terms in $w$'s by one while leaving the exponents
unchanged. The action by $g_2$ on the terms in $w_{2j}$'s, $j\geq 3$
(resp $w_{3j}$'s, $j\geq 4$) are given by reading the actions
$\n_1(w_{1,j-1})$'s (resp $\n_1(w_{2,j-1})$'s) on the formula $(5)$
(resp $(4)$) of Theorem $1$ of \cite{CL} and by incrementing all the
indices of the terms in $w$'s by one while leaving the exponents
unchanged. Next, we define

\begin{eqnarray}
\n_2(w_{12})&=&r\,w_{12}\\
\n_2(w_{13})&=&r\,w_{13}\\
\forall j\geq
4,\,\qquad\;\n_2(w_{1j})&=&\begin{split}m\,\wh{w_{3j}}+\ovs{\tp}{w_{3j}}-m\,r^{j-6}\,w_{12}
+m\,r^{j-5}\big(\unsurr\,\wh{w_{23}}-w_{13}\big)\\+m\,\big(\unsurr\,\wh{w_{2j}}-w_{1j}\big)\end{split}\\
\n_2(\ovs{\tp}{w_{3t}})&=&w_{1t}-m\,w_{2t}-m\,r\,w_{3t}+m\,r^{t-5}\,(w_{12}+r\,w_{13}-m\,r\,w_{23})\\&&\notag\\
\forall s\geq
4,\,\qquad\;\n_2(\ovs{\tp}{w_{st}})&=&\left\lb\begin{array}{l}
mr^{t-6}(r\wh{w_{2s}}+m\,w_{2s}-w_{1s})\\\\
+\;mr^{s-4}(r\wh{w_{2t}}+m\,w_{2t}-w_{1t})\\\\
+\;m\,r^{t-5}(\ovs{\tp}{w_{3s}}-r\,w_{3s})\\\\
+\;m\,r^{s-3}(\ovs{\tp}{w_{3t}}-r\,w_{3t})\\\\
+\;m^2(r^{s+t-10}+r^{s+t-8})(r\wh{w_{23}}+m\,w_{23}-w_{13}-\unsurr\,w_{12})\\\\
+\;r\,\ovs{\tp}{w_{st}}
\end{array}\right.\\
\forall s\geq
4,\,\qquad\;\n_2(\ovs{\tpp}{w_{st}})&=&\left\lb\begin{array}{l}
mr^{t-5}\bigg(\wh{w_{3s}}+r\,\ovs{\tp}{w_{3s}}+m(r+\unsurr)\,w_{3s}-(w_{1s}+r\,w_{2s})\bigg)\\\\
+\;mr^{s-3}\bigg(\wh{w_{3t}}+r\,\ovs{\tp}{w_{3t}}+m(r+\unsurr)\,w_{3t}-(w_{1t}+r\,w_{2t})\bigg)\\\\
+\;m^2(r^{s+t-10}+r^{s+t-8})\bigg(r^2(\wh{w_{23}}-w_{23})-(w_{12}+r\,w_{13})\bigg)\\\\
+\;r\,\ovs{\tpp}{w_{st}}
\end{array}\right.\\
\n_2(\ta)&=&\lr\\
\n_2(\lr)&=&\ta-m\,\lr\notag\\&&+\,\bigg[\frac{m}{lr}+\frac{(1-r^4)^2}{l}+\frac{m(1-r^4)(1-r^6)}{r^2}\bigg]\,\wh{w_{23}}
\end{eqnarray}

\noin When computing the action by $g_2$ on $w_{12}$ with the
tangles, one gets the same tangle as on the bottom line of Figure
$2$ (the tangle on the left hand side), except the top horizontal
strand over-crosses the vertical strand that it intersects. After
applying the commuting relation of the second type, the resulting
pole-related self-intersection has the opposite sign than the one on
Figure $2$. By Figure $1$, it yields a multiplication by $r$ instead
of a division by $r$. We thus obtain $(11)$.\\
From $(11)$, we easily derive $(12)$ by using $w_{13}=g_3\,w_{12}$
and by using the fact that $g_2$ and $g_3$ commute.\\
To get equality $(13)$, we use $w_{1j}=g_1^{-1}\,w_{2j}$. It follows
that
\begin{eqnarray*}
g_2w_{1j}&=&g_2g_1^{-1}\,w_{2j}\\
&=&g_1^{-1}g_2\,w_{2j}\\
&=&g_1^{-1}(\wh{w_{3j}}+m\,r^{j-5}(\wh{w_{23}}-w_{23})+m\,(\wh{w_{2j}}-w_{2j}))\\
&=&m\,\wh{w_{3j}}+\ovs{\tp}{w_{3j}}-m\,r^{j-6}\,w_{12}+m\,r^{j-5}\big(\unsurr\,\wh{w_{23}}-w_{13}\big)+m\,\big(\unsurr\,\wh{w_{2j}}-w_{1j}\big)
\end{eqnarray*}

\noin As for equality $(14)$, we have
\begin{eqnarray*}
g_2\ovs{\tp}{w_{3t}}&=&g_2g_1\wh{w_{3t}}\\
&=&g_1g_2\wh{w_{3t}}\\
&=&g_1(w_{2t}+m\,r^{t-4}\,w_{23}-m\,w_{3t})\\
&=&w_{1t}+m\,r^{t-3}\,w_{12}-m\,w_{2t}+m\,r^{t-4}(w_{13}+m\,w_{12}-m\,w_{23})-m\,r\,w_{3t}\\
&=&w_{1t}-m\,w_{2t}-m\,r\,w_{3t}+m\,r^{t-5}\,(w_{12}+r\,w_{13}-m\,r\,w_{23})
\end{eqnarray*}

\noin To derive Equality $(15)$, we use again
$g_2g_1\wh{w_{st}}=g_1g_2\wh{w_{st}}$. Then, we apply formula $(1)$
of \cite{CL} on $\n_1(\wh{w_{s-1,t-1}})$ and increment all the
indices of the result by one while leaving the exponents unchanged.
We thus get an expression for $g_2\,\wh{w_{st}}$. To finish, we
apply formulas $(3)$, $(4)$, $(5)$ of the current paper. \\
Applying $g_3$ on the result of $(15)$ now yields $(16)$. The action
by $g_3$ on $\ovs{\tp}{w_{3t}}$ appears later when $\n_3$ is defined
and is simply a multiplication by $r$ (see $(20)$).\\

\noin To compute $(18)$, we proceed as follows. First, we use the
Kauffman skein relation and we get
$$g_2\lr=\ta-m\,\lr+\frac{m}{l}\,e_2g_4g_3g_1\wh{w_{56}}$$
It remains to compute $e_2g_4g_3g_1\wh{w_{56}}$. We have
\begin{eqnarray*}
e_2g_4g_3g_1\wh{w_{56}}&=&e_2g_4g_2g_2^{-1}\ovs{\tpp}{w_{56}}\\
&=&e_2e_4g_2^{-1}\ovs{\tpp}{w_{56}}
\end{eqnarray*}

\noin To compute $g_2^{-1}\ovs{\tpp}{w_{56}}$, we use the formulas
$(11),(12),(13),(16)$ above and the formulas $(2)-(7)$ of \cite{CL}.
The computations are long. We content ourselves to give the result
but show on an example how they work. In the following equalities,
the symbol $\blacklozenge$ denotes the action for the representation
defined in Theorem $1$ of \cite{CL} and the symbol $<<>>$ is the
operation that increments the indices by one while leaving the
exponents unchanged. We will also use these notations whenever we
need them later along the paper. So for instance, we have

\begin{eqnarray*}
g_2^{-1}w_{15}&=&g_2^{-1}\,g_1^{-1}\,w_{25}\\
&=&g_1^{-1}g_2^{-1}\,w_{25}\\
&=&g_1^{-1}\,<<g_1^{-1}\bl\,w_{14}>>\\
&=&g_1^{-1}\,<<\wh{w_{24}}-m\,w_{12}+m\,\wh{w_{14}}>>\\
&=&g_1^{-1}\,(\wh{w_{35}}-m\,w_{23}+m\,\wh{w_{25}})\\
&=&\ovs{\tp}{w_{35}}+m\,\wh{w_{35}}-\frac{m}{r}\,w_{12}-m\,w_{13}+\frac{m}{r}\,\wh{w_{25}}
\end{eqnarray*}

\noin By doing such computations, we get

\begin{equation}
g_2^{-1}\ovs{\tpp}{w_{56}}=\begin{split}\unsurr\bigg[\ovs{\tpp}{w_{56}}-m^2(r+r^3)(r^2\,l\,\wh{w_{23}}-r\,w_{23}-\unsurr\,w_{12}-w_{13})\\
-m\,r^2\left(\begin{array}{l}r\,w_{16}+r^2\,w_{26}-m\,(1+r^2)\,w_{36}-r\,\wh{w_{36}}\\
-m\,l\,r^3(1+r^2)\,\wh{w_{23}}+m\,r^2(1+r^2)\,w_{23}\\+m\,(1+r^2)\,w_{12}+m\,r(1+r^2)\,w_{13}-r^2\,\ovs{\tp}{w_{36}}
\end{array}\right)\\
-m\,r\left(\begin{array}{l}r\,w_{15}+r^2\,w_{25}-m\,(1+r^2)\,w_{35}-r\,\wh{w_{35}}\\
-m\,l\,r^2(1+r^2)\,\wh{w_{23}}+(1-r^4)\,w_{23}\\+m\,(\unsurr+r)\,w_{12}+m(1+r^2)\,w_{13}-r^2\,\ovs{\tp}{w_{35}}
\end{array}\right) \bigg]\end{split}
\end{equation}

\noin Next, we must act on this result with $e_2e_4$. The actions
are summarized below.
\newtheorem{Lemma}{Lemma}
\begin{Lemma} The following equalities hold.
$$\begin{array}{ccc}
\begin{array}{ccc}
e_2e_4\,w_{16}&=&0\\
e_2e_4\,w_{26}&=&0\\
e_2e_4\,w_{36}&=&l\,r\,\wh{w_{23}}\\
e_2e_4\,\wh{w_{36}}&=&\wh{w_{23}}\\
e_2e_4\,\ovs{\tp}{w_{36}}&=&r\,\wh{w_{23}}\\
e_2e_4\,\ovs{\tpp}{w_{56}}&=&\wh{w_{23}}
\end{array}&\begin{array}{ccc}
e_2e_4\,w_{15}&=&0\\
e_2e_4\,w_{25}&=&0\\
e_2e_4\,w_{35}&=&l\,\wh{w_{23}}\\
e_2e_4\,\wh{w_{35}}&=&\unsurr\,\wh{w_{23}}\\
e_2e_4\,\ovs{\tp}{w_{35}}&=&\wh{w_{23}}
\end{array}
&\begin{array}{ccc}
e_2e_4\,w_{12}&=&0\\
e_2e_4\,w_{13}&=&\unsurr\,\wh{w_{23}}\\
e_2e_4\,w_{23}&=&\wh{w_{23}}\\
e_2e_4\,\wh{w_{23}}&=&\wh{w_{23}}\\
\end{array}
\end{array}$$
\end{Lemma}

\noin\textsc{Proof.} In most cases, we use the tangles with
\begin{center}
\epsfig{file=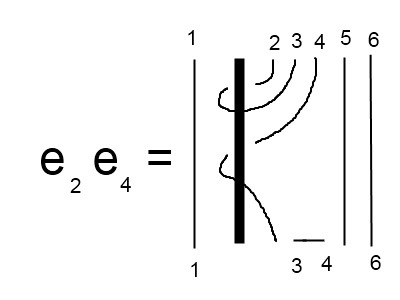, height=4.5cm}
\end{center}
In some cases it is easier to proceed algebraically like below.
\begin{eqnarray*}
e_2e_4\,\ovs{\tp}{w_{36}}&=&e_2e_4g_1\wh{w_{36}}\\
&=&g_1e_2e_4\wh{w_{36}}\\
&=&g_1\,\wh{w_{23}}\\
&=&r\,\wh{w_{23}}
\end{eqnarray*}
By $(A1)$ and $(P)$, the computation reduces to computing
$e_2e_4\wh{w_{36}}$, which is easily computed with the tangles and
is simply $\wh{w_{23}}$. The last equality then holds by $(4)$. Of
course this equality can also be visualized by using the tangles.\\
Also,
\begin{eqnarray*}e_2e_4w_{12}&=&e_2(e_4e_1)e_3e_4e_5e_6\\
&=& 0\;\; \text{in $M$}\end{eqnarray*}

\noin A more tricky computation is the following.
\begin{eqnarray*}
e_2e_4\,\ovs{\tpp}{w_{56}}&=&e_2e_4g_3g_1\wh{w_{56}}\\
&=&e_2e_4g_3g_4g_4^{-1}g_1\wh{w_{56}}\\
&=&e_2e_4e_3g_1g_4^{-1}\wh{w_{56}}\\
&=&\unsurr\,e_2e_4e_3\ovs{\tp}{w_{56}}
\end{eqnarray*}
Then, by using
$$e_3=1+\unsur{m}(g_3-g_3^{-1}),$$
and forthcoming equalities $(21)$ and $(22)$, we get
$$e_3\,\ovs{\tp}{w_{56}}=r\,w_{23}$$

\noin And since $e_2e_4\,w_{23}=\wh{w_{23}}$, we get
$e_2e_4\,\ovs{\tpp}{w_{56}}=\wh{w_{23}}$, as in the lemma. We are
now ready to apply $e_2e_4$ on the right hand side of $(19)$. We get

$$e_2e_4g_2^{-1}\,\ovs{\tpp}{w_{56}}=\bigg[\unsurr+\bigg(\unsur{r^2}-r^2\bigg)\bigg(l(1-r^6)+r^3+r^5\bigg)\bigg]\wh{w_{23}}$$

\noin It remains to multiply by $\frac{m}{l}$ and it yields the last term of $(18)$. \\



\noin We now move on to $\n_3$. The action by $g_3$ on the terms
$w_{ij}$'s with $i\geq 2$ is the same action as in the
representation of \cite{CL}. The action by $g_3$ on the terms
$\wh{w_{ij}}$'s is obtained as follows. Read the action of $g_2$ on
the term $\wh{w_{i-1,j-1}}$ using formulas $(8)-(10)$, $(12)$,
$(13)$, $(16)$ of \cite{CL}, then increment all the indices of the
result by one without modifying the exponents. The other actions are
given by

\begin{eqnarray}
\n_3(\ovs{\tp}{w_{3t}})&=&r\,\ovs{\tp}{w_{3t}}\\
\forall s\geq 4,\,\qquad\;\n_3(\ovs{\tp}{w_{st}})&=&\ovs{\tpp}{w_{st}}\\
\forall s\geq 4,\,\qquad\;\n_3(\ovs{\tpp}{w_{st}})&=&\ovs{\tp}{w_{st}}-m\,\ovs{\tpp}{w_{st}}+\frac{m}{l}\,r^{s+t-10}\,w_{23}\\
\n_3(\ta)&=&r\,\ta\\
\n_3(\lr)&=&r\,\lr
\end{eqnarray}

\noin To get equality $(20)$, we notice that

\begin{eqnarray*}
g_3\ovs{\tp}{w_{3t}}&=&g_3g_1g_3\wh{w_{2t}}\\
&=&g_1g_3g_1\wh{w_{2t}}\\
&=&r\,g_1\wh{w_{3t}}\\
&=&r\,\ovs{\tp}{w_{3t}}
\end{eqnarray*}

\noin To get Equality $(22)$,
apply the Kauffman skein relation and get the first two terms in
$(22)$, plus a multiple of $e_3g_1\wh{w_{st}}$. We compute the
latter term as follows. We have
\begin{eqnarray*}
e_3g_1\wh{w_{st}}&=&e_3g_1g_3g_3^{-1}\,\wh{w_{st}}\\
&=&\unsurr\,e_3e_1\wh{w_{st}}\\
&=&\unsurr\,e_3(r^{s+t-9}\,w_{12})\\
&=&r^{s+t-10}\,w_{23}
\end{eqnarray*}


\noin Equality $(24)$ is obtained by using the sequence of relations
\begin{eqnarray*}
g_3\lr&=&g_3g_2g_4g_3g_1\wh{w_{56}}\\
&=&g_2g_3g_4g_3g_1\wh{w_{56}}\\
&=&g_2g_4g_3g_4g_1\wh{w_{56}}\\
&=&g_2g_4g_3g_1g_4\wh{w_{56}}\\
&=&r\,\lr
\end{eqnarray*}

\noin We go on with $\n_4$. We only provide below the actions that
don't easily follow from the $D_5$ case in \cite{CL} or which are
the most difficult.

\begin{eqnarray}
\forall\,t\in\lb 5,6\rb,\qquad\n_4(\ovs{\tp}{w_{4t}})&=&\ovs{\tp}{w_{3t}}+\frac{m\,r^{t-5}}{l}\,w_{34}-m\,\ovs{\tp}{w_{4t}}\\
\forall\,t\in\lb 5,6\rb,\qquad\n_4(\ovs{\tpp}{w_{4t}})&=&r\,\ovs{\tpp}{w_{4t}}\\
\n_4(\ta)&=&\ovs{\tpp}{w_{56}}-m\,\ta+\frac{m}{l}\,w_{34}\\&&\notag\\
\n_4(\lr)&=&\left\lb\begin{array}{l}m^2(r+r^3)\Bigg(\begin{split}w_{45}-r\,w_{35}+r(w_{46}-r\,w_{36})\\-r^2[w_{14}-r\,w_{13}+r(w_{24}-r\,w_{23})]\\-w_{34}\end{split}\;\Bigg)\\\\\\
+\,m\,r^2\big(\wh{w_{45}}-r\,\wh{w_{35}}+r\,(\ovs{\tp}{w_{45}}-r\,\ovs{\tp}{w_{35}}+\wh{w_{46}}-r\,\wh{w_{36}})+r^2(\ovs{\tp}{w_{46}}-r\,\ovs{\tp}{w_{36}})\big)\\\\
+\,m^2\,r^2(1+r^2)^2\,(\wh{w_{24}}-r\,\wh{w_{23}})\\+\,r\,\lr\end{array}
\right.\end{eqnarray}

\noin We now comment on these equalities. Equality $(25)$ is done by
using the fact that $g_1$ and $g_4$
commute. \\
Equality $(26)$ is derived by using a trick already used before.
\begin{eqnarray*}g_4\ovs{\tpp}{w_{4t}}&=&g_4g_3g_4g_4^{-1}g_1\,\wh{w_{4t}}\\
&=&g_3g_4g_3g_1g_4^{-1}\,\wh{w_{4t}}\\
&=&g_3g_4g_3\ovs{\tp}{w_{3t}}\\
&=&r\,g_3g_4\ovs{\tp}{w_{3t}}\\
&=&r\,g_3\ovs{\tp}{w_{4t}}\\
&=&r\,\ovs{\tpp}{w_{4t}}
\end{eqnarray*}
\noin For the last term of $(27)$, we use the mixed Artin relation
$(MA)$ like we did when we computed $\n_3(\ovs{\tpp}{w_{st}})$. This
computation is left to the reader. Finally, we have

\begin{eqnarray*}
g_4.\,\lr&=&g_4g_2g_4.\,\ovs{\tpp}{w_{56}}\\\\
&=&g_2g_4g_2.\,\ovs{\tpp}{w_{56}}\\\\
&=&\begin{split}g_2.\,\bigg\lb m\,r\,\bigg(\wh{w_{45}}+r\,\ovs{\tp}{w_{45}}+m(\unsurr+r)\,w_{45}-r\,w_{15}-r^2\,w_{25}\bigg)\\
+mr^2\,\bigg(\wh{w_{46}}+r\,\ovs{\tp}{w_{46}}+m(r+\unsurr)\,w_{46}-r\,w_{16}-r^2\,w_{26}\bigg)\\
+m^2(r+r^3)\,\bigg(r^2(\wh{w_{24}}-w_{24})-r\,w_{12}-r\,w_{14}\bigg)\\
+\,r\,\ta\bigg\rb\end{split}
\end{eqnarray*}

\noin To compute the term in the last equality, we used the matrix
of the left action by $g_2$ and Maple. We obtained $(28)$.\\

\noin Again, for $\n_5$, we only provide here the actions that are
the less straightforward. The following equalities hold.

\begin{eqnarray}
\n_5(\ovs{\tp}{w_{35}})&=&\ovs{\tp}{w_{34}}+\frac{m\,r}{l}\,w_{45}-m\,\ovs{\tp}{w_{35}}\\
\n_5(\ovs{\tp}{w_{56}})&=&\ovs{\tp}{w_{46}}+\frac{m}{l}\,w_{45}-m\,\ovs{\tp}{w_{56}}\\
\n_5(\ovs{\tpp}{w_{45}})&=&r\,\ovs{\tpp}{w_{45}}\\
\n_5(\ovs{\tpp}{w_{46}})&=&\ovs{\tpp}{w_{56}}\\
\n_5(\ovs{\tpp}{w_{56}})&=&\ovs{\tpp}{w_{46}}+\frac{mr}{l}\,w_{45}-m\,\ovs{\tpp}{w_{56}}\\
\n_5(\ta)&=&r\,\ta\\
\n_5(\lr)&=&r\,\lr
\end{eqnarray}
\noin The following sets of equalities explain $(33)$, $(34)$ and
$(35)$ respectively.
\begin{eqnarray*}
\n_5(\ovs{\tpp}{w_{56}})&=&g_5g_3\ovs{\tp}{w_{56}}\\
&=&g_3g_5\ovs{\tp}{w_{56}}\\
&=&g_3(\ovs{\tp}{w_{46}}+\frac{m}{l}\,w_{45}-m\,\ovs{\tp}{w_{56}})\\
&=&\ovs{\tpp}{w_{46}}+\frac{mr}{l}\,w_{45}-m\,\ovs{\tpp}{w_{56}}
\end{eqnarray*}
\begin{eqnarray*}
g_5(\ta)&=&g_5g_4\ovs{\tpp}{w_{56}}\\
&=&g_5g_4g_5g_5^{-1}\ovs{\tpp}{w_{56}}\\
&=&g_4g_5g_4\ovs{\tpp}{w_{46}}\\
&=&g_4g_5(r\,\ovs{\tpp}{w_{46}})\\
&=&r\,g_4\ovs{\tpp}{w_{56}}\\
&=&r\,\ta
\end{eqnarray*}
\begin{eqnarray*}
g_5(\lr)&=&g_5g_2\ta\\
&=&g_2g_5\ta\\
&=&r\,g_2\ta\\
&=&r\,\lr
\end{eqnarray*}

\noin The action by $g_6$ is in many ways similar to the one by
$g_5$ and is left to the reader. We use in particular the formulas
$(13)$ and $(15)$ of \cite{CL} with $j=6$.\\\\
With Maple, we formed the matrices of the $\n_i$'s and checked that
they verify all the Artin relations. Further, we define
$\n(e_i)=\frac{l}{m}(\n_i^2+m\,\n_i-id_{V_6})$ and we check that the
matrices of the $\n(e_i)$'s and the matrices of the $\n_i$'s satisfy
the delooping relations $(DL1)$ and $(DL2)$. We thus obtain a
representation of the BMW algebra of type $E_6$. In the following
part, we study its reducibility.

\section{Reducibility}
The reducibility result is based on the fundamental following lemma.
\begin{Lemma}
Let $\U$ be a proper invariant subspace of $V_6$. Then $\U$ is
annihilated by all the $e_i$'s. And so, it is also annihilated by
all the conjugates of the $e_i$'s.
\end{Lemma}
\noin \textsc{Proof.} An action by $e_1$ (resp $e_2$, resp $e_i$,
$i\geq 3$) on any basis vector always yields a multiple of $w_{12}$
(resp $\wh{w_{23}}$, resp $w_{i-1,i}$). Moreover, if one of these
vectors lies in $\U$, then all the other basis vectors are also in
$\U$ and $\U$ would hence be the whole space. This contradicts $\U$
is proper.

\begin{Definition}
Define 
$$\begin{array}{cccc}
X_{13}=g_3\,e_1\,g_3^{-1}&\text{and}&X_{1j}=g_j\,X_{1,j-1}\,g_j^{-1}&\forall\,4\leq
j\leq 6\\
X_{24}=g_4\,e_3\,g_4^{-1}&\text{and}&X_{2j}=g_j\,X_{2,j-1}\,g_j^{-1}&\forall\,j\in\lb
5,6\rb\\
X_{35}=g_5\,e_4\,g_5^{-1}&\text{and}&X_{36}=g_6\,X_{35}\,g_6^{-1}&\\
X_{46}=g_6\,e_5\,g_6^{-1}&&&
\end{array}$$

Define also $$\begin{array}{cccc} 
\wh{X_{24}}=g_4\,e_2\,g_4^{-1}&\text{and}&\wh{X_{2j}}=g_j\,\wh{X_{2,j-1}}\,g_j^{-1}&\forall
4\leq j\leq 6\\
\wh{X_{3j}}=g_3\,\wh{X_{2j}}\,g_3^{-1}&\forall 4\leq j\leq 6&\\
\wh{X_{4j}}=g_4\,\wh{X_{3j}}\,g_4^{-1}&\forall\,j\in\lb 5,6\rb&\\
\wh{X_{56}}=g_5\,\wh{X_{46}}\,g_5^{-1}
\end{array}$$

Define further $$\begin{array}{cc}
\ovs{\tp}{X_{st}}=g_1\,\wh{X_{st}}\,g_1^{-1}&\forall\,s\geq 3\\
\ovs{\tpp}{X_{st}}=g_3\,\ovs{\tp}{X_{st}}\,g_3^{-1}&\forall\,s\geq 4\\
\ovs{\tppp}{X_{56}}=g_4\,\ovs{\tpp}{X_{56}}\,g_4^{-1}&\\
\ovs{\triangle^{l}}{X_{56}}=g_2\,\ovs{\tppp}{X_{56}}\,g_2^{-1}&
\end{array}\\$$

Define finally $$\begin{array}{l}
X_{12}=e_1\;\;\text{and}\;\;X_{23}=e_3\;\;\text{and}\;\;X_{34}=e_4\;\;\text{and}\;\;X_{45}=e_5\;\;\text{and}\;\;X_{56}=e_6\;\;\text{and}\;\;
\wh{X_{23}}=e_2\end{array}$$

Let $S$ denote the sum of all these conjugates.
\end{Definition}

\noin Suppose $\n^{(6)}$ is reducible and let $\U$ be a proper
invariant subspace of $V_6$. By Lemma $2$, we have
$$S.\,\U=0$$
\noin Then, the determinant of the left action by $S$ must be zero,
as $\U$ is non-trivial. We computed the determinant of this action
with Mathematica and we obtained
$$det(S)=\frac{(l-r^3)^{15}(l+r^3)^{30}(-1+l\,r^3)^{20}(1+l\,r^9)^6(-1+l\,r^{21})}{l^{36}\,r^{99}(-1+r^2)^{36}}$$
\noin This yields a necessary condition for the representation
$\n^{(6)}$ to be reducible. For the converse, we show the following
result.
\begin{Proposition}
The $F$-vector space
$$\mathcal{K}_6=\begin{array}{l}\bigcap_{1\leq i<j\leq
6}\text{Ker}\,\n^{(6)}(X_{ij})\;\cap\;\bigcap_{2\leq i<j\leq
6}\text{Ker}\,\n^{(6)}(\wh{X_{ij}})\;\cap\;\bigcap_{3\leq i<j\leq
6}\text{Ker}\,\n^{(6)}(\ovs{\tp}{X_{ij}})\\\\\cap\;\bigcap_{4\leq
i<j\leq
6}\text{Ker}\,\n^{(6)}(\ovs{\tpp}{X_{ij}})\;\cap\;\text{Ker}\,\n^{(6)}(\ovs{\tppp}{X_{56}})\;
\cap\text{Ker}\,\n^{(6)}(\ovs{\triangle^{l}}{X_{56}})\end{array}$$
is a $BMW(E_6)$-submodule of $V_6$.
\end{Proposition}

\noin \textsc{Proof.} Given any vector $x$ in $\ab$, by the work
already achieved in the proofs of Lemma $5$ of \cite{CL} and of
Lemma $2$ in Chapter $2$ of \cite{THE}, it suffices to show that
$g_1\,x$ is annihilated by the $X_{st}$ terms which either carry a
hat or a triangle, that $g_2\,x$ is annihilated by all the
$X_{1t}$'s with $t\geq 2$ and that when $i\geq 2$, $g_i\,x$ is
annihilated
by all the $X_{st}$ terms which wear a triangle. 
First, since $g_1$ commutes with the $\wh{X_{2j}}$'s, we have
$g_1\,x\in\cap_{3\leq j\leq 6}\text{Ker}\,\n^{(6)}(\wh{X_{2j}})$. To
show that $g_1\,x$ is in the kernel of $\wh{X_{s,j}}$ with $s\geq
3$, it suffices to show that $g_1^{-1}\,x$ is in the kernel of
$\wh{X_{s,j}}$. Notice
$$\wh{X_{s,j}}\,g_1^{-1}\,x=g_1^{-1}\ovs{\tp}{X_{s,j}}\,x=0$$ So,
\begin{equation}g_1\,x\in\cap_{2\leq i<j\leq
6}\,\text{Ker}\,\n^{(6)}(\wh{X_{ij}})\end{equation}\noin Next, by
definition of $\ovs{\tp}{X_{st}}$, we have
$g_1\,x\in\text{Ker}\,\n^{(6)}(\ovs{\tp}{X_{st}})$ for all $s\geq
3$. Fix $s\geq 4$. We now show that $g_1^{-1}$ commutes with
$\ovs{\tpp}{X_{st}}$. We have
\begin{eqnarray}
\ovs{\tpp}{X_{st}}\,g_1^{-1}\,x&=&g_1^{-1}\,g_1\,g_3\,\ovs{\tp}{X_{st}}\,\,g_3^{-1}g_1^{-1}\,x\\
&=&g_1^{-1}g_1\,g_3\,g_1\,\wh{X_{st}}\,g_1^{-1}g_3^{-1}g_1^{-1}\,x\\
&=&g_1^{-1}g_3\,g_1\,g_3\,\wh{X_{st}}\,g_3^{-1}g_1^{-1}g_3^{-1}\,x\\
&=&g_1^{-1}g_3\,g_1\,\wh{X_{st}}\,g_1^{-1}g_3^{-1}\,x\\
&=&g_1^{-1}\ovs{\tpp}{X_{st}}\,x\\
&=&0
\end{eqnarray}
Equality $(40)$ is easy to see on the tangles after doing a
Reidemeister's move of type \setcounter{co}{3} (\Roman{co}).
Further,
\begin{eqnarray}
\ovs{\tppp}{X_{56}}\,g_1^{-1}\,x&=&g_4\,g_3\,g_1\,\wh{X_{56}}\,\,g_1^{-1}g_3^{-1}g_1^{-1}g_4^{-1}\,x\\
&=&g_1^{-1}g_4\,g_3\,g_1\,g_3\,\wh{X_{56}}\,g_3^{-1}g_1^{-1}g_3^{-1}g_4^{-1}\,x\\
&=&g_1^{-1}g_4\,g_3\,g_1\,\wh{X_{56}}\,g_1^{-1}g_3^{-1}g_4^{-1}\,x\\
&=&g_1^{-1}\ovs{\tppp}{X_{56}}\,x\\
&=&0
\end{eqnarray}
\noin Similarly, we show that
$$\ovs{\triangle^{l}}{X_{56}}\,\,g_1^{-1}\,x=g_1^{-1}\,\ovs{\triangle^l}{X_{56}}\,x=0$$
\noin At this stage, we conclude that $g_1\,x\in\ab$. Next, we show
that $g_2\,x$ is annihilated by all the $X_{1t}$'s with $t\geq 2$.
We have $X_{12}\,g_2\,x=0$ since $e_1$ commutes with $g_2$. Further,
when $t\geq 3$, we have $X_{1t}=g_1^{-1}\,X_{2t}\,g_1$, so that
$$X_{1t}\,g_2\,x=g_1^{-1}\,X_{2t}\,g_2\,(g_1\,x)$$
We just showed above that $g_1\,x\in\ab$. It follows that
$g_2(g_1\,x)\in\text{Ker}(X_{2t})$ by Lemma $5$ of \cite{CL}.

\noin It remains to show that when $i\geq 2$, we have $g_i\,x$
belongs to the kernel of the $\ovs{\triangle}{X_{st}}$'s. We begin
with $i=2$. Obviously,
$g_2\,x\in\text{Ker}\,\n^{(6)}(\ovs{\triangle^l}{X_{56}})$ and
$g_2^{-1}\,x\in\text{Ker}\,\n^{(6)}(\ovs{\tppp}{X_{56}})$.\\
Next, for every $s\geq 3$, we have
$$\ovs{\tp}{X_{st}}\,g_2^{-1}\,x=g_1\,\wh{X_{st}}\,g_2^{-1}\,g_1^{-1}\,x$$
As we have seen above that $g_1^{-1}\,x\in\ab$, the member to the
right hand side of the latter equality is zero. Before we can
proceed the other $\ovs{\triangle}{X}$ terms, we must deal with the
case $i=3$ first.
\noin We use the following equalities, with the same techniques as
before.
\begin{Lemma}
\begin{eqnarray*}
\hspace{-2cm}\ovs{\tp}{X_{st}}\,\,g_3^{-1}\,x&=&\begin{cases}g_3^{-1}\,\ovs{\tpp}{X_{st}}\,x&\text{if $s\geq 4$}\\
g_1\,g_3\,g_1^{-1}\,\wh{X_{2t}}\,g_3^{-1}g_1^{-1}\,x&\text{if $s=3$}
\end{cases}\\
\hspace{-2cm}\ovs{\tpp}{X_{st}}\,g_3\,x&=&g_3\,\,\ovs{\tp}{X_{st}}\,x\\
\hspace{-2cm}\ovs{\tppp}{X_{56}}\,g_3^{-1}\,x&=&g_3^{-1}\,\ovs{\tppp}{X_{56}}\,x\\
\hspace{-2cm}\ovs{\triangle^l}{X_{56}}\,\,g_3^{-1}\,x&=&g_3^{-1}\,\ovs{\triangle^l}{X_{56}}\,x
\end{eqnarray*}
\end{Lemma}
\noin The second part of the first equality is obtained as follows.
\begin{eqnarray*}
\ovs{\tp}{X_{3t}}\,\,g_3^{-1}\,x&=&g_1\,g_3\,\wh{X_{2t}}\,g_3^{-1}g_1^{-1}g_3^{-1}\,x\\
&=&g_1\,g_3\,\wh{X_{2t}}\,g_1^{-1}g_3^{-1}g_1^{-1}\,x\\
&=&g_1\,g_3\,g_1^{-1}\,\wh{X_{2t}}\,g_3^{-1}g_1^{-1}\,x
\end{eqnarray*}
Because we have seen earlier that $g_1^{-1}\,x\in\ab$, we know that
$g_3^{-1}(g_1^{-1}\,x)\in\,\text{Ker}(\wh{X_{2t}})$. Thus, all the
results above are zero. This ends the case $i=3$. It is now easy to
close the case $i=2$. Indeed, we have for every $s\geq 4$,
$$\ovs{\tpp}{X_{st}}\,g_2^{-1}\,x=g_3\,g_1\,\wh{X_{st}}\,g_2^{-1}\,g_1^{-1}g_3^{-1}\,x$$
By the case $i=3$, we have $g_3^{-1}\,x\in\ab$. Then by the case
$i=1$, we get $y=g_1^{-1}g_3^{-1}\,x\in\ab$. It follows that
$g_2^{-1}\,y\in\text{Ker}(\wh{X_{st}})$.\\
We now deal with the case $i=4$. When $i=4$, the following
equalities hold.
\begin{Lemma}
\begin{eqnarray}\hspace{-2cm}
\ovs{\tpp}{X_{st}}\,g_4^{-1}\,x&=&\begin{cases}g_4^{-1}\,\ovs{\tppp}{X_{56}}\,x&\text{if $s=5$ and $t=6$}\notag\\
g_3\,g_4\,g_1\,g_3\,g_1^{-1}\,\wh{X_{2t}}\,g_3^{-1}g_4^{-1}g_1^{-1}g_3^{-1}\,x&\text{if $s=4$}\end{cases}\notag\\
\hspace{-2cm}\ovs{\tppp}{X_{56}}\,g_4\,x&=&g_4\,\ovs{\tpp}{X_{56}}\,x\notag\\
\hspace{-2cm}\ovs{\triangle^l}{X_{56}}\,g_4^{-1}\,x&=&g_2\,g_4\,g_3\,g_1\,\wh{X_{56}}\,g_2^{-1}g_1^{-1}g_3^{-1}g_4^{-1}g_2^{-1}\,x\\
\hspace{-2cm}&&\notag\\
\hspace{-2cm}\ovs{\tp}{X_{st}}\,g_4^{-1}\,x&=&g_4^{-1}g_1\,g_4\,\wh{X_{st}}\,g_4^{-1}g_1^{-1}\,x=
\begin{cases}
g_4^{-1}\,\ovs{\tp}{X_{34}}\,x&\text{if $s=3$ and $t=4$}\\
g_4^{-1}\,\ovs{\tp}{X_{4t}}\,x&\text{if $s=3$ and $t\geq 5$}\\
g_4^{-1}\,\ovs{\tp}{X_{56}}\,\,x&\text{if $s=5$ and
$t=6$}\end{cases}\notag\\
\hspace{-2cm}\ovs{\tp}{X_{4t}}\,g_4\,x&=&g_4\,\ovs{\tp}{X_{3t}}\,x\qquad\forall\,t\in\lb
5,6\rb\notag
\end{eqnarray}
\end{Lemma}
\noin To get the second term in the first equality, proceed as
follows.
\begin{eqnarray*}
\ovs{\tpp}{X_{4t}}\,g_4^{-1}\,x&=&g_3\,g_1\,g_4\,\wh{X_{3t}}\,g_4^{-1}g_1^{-1}g_3^{-1}g_4^{-1}\,x\\
&=&g_3\,g_4\,g_1\,g_3\,\wh{X_{2t}}\,g_3^{-1}g_1^{-1}g_3^{-1}g_4^{-1}g_3^{-1}\,x\\
&=&g_3\,g_4\,g_1\,g_3\,\wh{X_{2t}}\,g_1^{-1}g_3^{-1}g_1^{-1}g_4^{-1}g_3^{-1}\,x\\
&=&g_3\,g_4\,g_1\,g_3\,g_1^{-1}\,\wh{X_{2t}}\,g_3^{-1}g_4^{-1}g_1^{-1}g_3^{-1}\,x
\end{eqnarray*}
\noin It is now easy to conclude. By Lemma $3$, we know that
$g_3^{-1}\,x\in\ab$. Then, $y=g_1^{-1}(g_3^{-1}\,x)$ is also in
$\ab$. In particular, $y$ belongs to the kernel of $X_{ij}$ and of
$\wh{X_{ij}}$ for all the integers $i$ and $j$ with $2\leq i<j\leq
6$. Then, an adequate adaptation of Lemma $5$ of \cite{CL} shows
that
$$g_3^{-1}g_4^{-1}\,y\in\,\text{Ker}\,(\wh{X_{2t}})$$
\noin To show that the right member of equality $(48)$ is zero, it
suffices to show that
$$g_1^{-1}g_3^{-1}g_4^{-1}g_2^{-1}\,x\in\cap_{2\leq i<j\leq
6}\text{Ker}\,(X_{ij}\cap\wh{X_{ij}})$$ For convenience, denote by
$\ovl{\mathcal{K}_5}$ this intersection. First, we have
$g_2^{-1}\,x\in\ab$. By the equalities of Lemma $4$, it follows that
$g_4^{-1}g_2^{-1}\,x$ belongs to the kernel of all the $X$'s, except
possibly $\ovs{\triangle^l}{X_{56}}$. Then,
$g_1^{-1}(g_4^{-1}g_2^{-1}\,x)\in\ovl{\mathcal{K}_5}$. It follows
that $g_3^{-1}g_4^{-1}g_2^{-1}\,x$ is in the kernel of all the
$X$'s, except possibly $\ovs{\triangle^l}{X_{56}}$. Then,
$g_1^{-1}g_3^{-1}g_4^{-1}g_2^{-1}\,x$ belongs to
$\ovl{\mathcal{K}_5}$ as desired. This ends the case $i=4$.

\begin{Lemma}
\begin{eqnarray}
\hspace{-2cm}\ovs{\tp}{X_{35}}\,g_5\,x&=&g_5\,\ovs{\tp}{X_{34}}\,x\notag\\
\hspace{-2cm}\ovs{\tp}{X_{56}}\,g_5\,x&=&g_5\,\ovs{\tp}{X_{46}}\,x\notag\\
\hspace{-2cm}\ovs{\tp}{X_{34}}\,g_5^{-1}\,x&=&g_5^{-1}\,\ovs{\tp}{X_{35}}\,x\notag\\
\hspace{-2cm}\ovs{\tp}{X_{46}}\,g_5^{-1}\,x&=&g_5^{-1}\,\ovs{\tp}{X_{56}}\,x\notag\\
\hspace{-2cm}\ovs{\tp}{X_{st}}\,g_5^{-1}\,x&=&g_5^{-1}\,\ovs{\tp}{X_{st}}\,x\qquad\text{if
$(s,t)\in\lb (4,5),(3,6)\rb$}\notag\\
\hspace{-2cm}\ovs{\tpp}{X_{56}}\,g_5\,x&=&g_5\,\ovs{\tpp}{X_{46}}\,x\notag\\
\hspace{-2cm}\ovs{\tpp}{X_{45}}\,g_5\,x&=&g_5\,\ovs{\tpp}{X_{45}}\,x\notag\\
\hspace{-2cm}\ovs{\tpp}{X_{46}}\,g_5^{-1}\,x&=&g_5^{-1}\,\ovs{\tpp}{X_{56}}\,x\notag\\
\hspace{-2cm}\ovs{\tppp}{X_{56}}\,g_5^{-1}\,x&=&g_4\,g_5\,g_3\,g_4\,g_1\,g_3\,g_1^{-1}\,\wh{X_{26}}\,g_3^{-1}g_4^{-1}g_5^{-1}g_1^{-1}g_3^{-1}g_4^{-1}\,x\\
\hspace{-2cm}\ovs{\triangle^l}{X_{56}}\,g_5^{-1}\,x&=&g_2\,\ovs{\tppp}{X_{56}}\,g_5^{-1}g_2^{-1}\,x
\end{eqnarray}
\end{Lemma}
\vspace{0.05cm} \noin Equality $(49)$ is obtained by using the
expression for $\ovs{\tpp}{X_{46}}\,g_4^{-1}$ from Lemma $4$. \\With
$(50)$, the fact that $\ovs{\triangle^l}{X_{56}}\,g_5^{-1}\,x$ is
zero follows from $(49)$ and from the fact that $g_2^{-1}\,x\in\ab$.
It remains to deal with $i=6$.

\begin{Lemma}
\begin{eqnarray*}
\ovs{\triangle^l}{X_{56}}\,\,\,g_6\,x&=&g_6\,\ovs{\triangle^l}{X_{56}}\,x\\
\ovs{\tppp}{X_{56}}\,g_6\,x&=&g_6\,\ovs{\tppp}{X_{56}}\,x\\
\hspace{-2cm}\ovs{\tpp}{X_{45}}\,\,g_6^{-1}\,x&=&g_6^{-1}\,\ovs{\tpp}{X_{46}}\,x\\
\ovs{\tpp}{X_{s,t}}\,g_6\,x&=&
\begin{cases}g_6\,\ovs{\tpp}{X_{56}}\,x
&\text{if $s=5$ and
$t=6$}\\
g_6\,\ovs{\tpp}{X_{45}}\,x&\text{if $s=4$ and $t=6$}\\
\end{cases}\\
\hspace{-2cm}\ovs{\tp}{X_{s,5}}\,\,g_6^{-1}\,x&=&g_6^{-1}\,\ovs{\tp}{X_{s,6}}\,x\qquad\forall\,s\in\lb
3,4\rb\\
\hspace{-2cm}\ovs{\tp}{X_{s,6}}\,\,g_6\,x&=&g_6\,\,\ovs{\tp}{X_{s,5}}\,x\qquad\forall\,s\in\lb
3,4\rb\\
\hspace{-2cm}\ovs{\tp}{X_{st}}\,\,g_6\,x&=&g_6\,\ovs{\tp}{X_{st}}\,x\qquad\text{if
$(s,t)\in\lb(3,4),(5,6)\rb$}\\
\end{eqnarray*}
\end{Lemma}
\noin This achieves the proof of Proposition $3$. A consequence of
Lemma $2$ and Proposition $3$ is that when the representation is
reducible, any proper invariant subspace of $V_6$ is a
$BMW(E_6)$-submodule of $\ab$. The following proposition studies the
dimension of $\ab$ as a vector space over $F$.


\begin{Proposition}
$\text{Ker}\,S\subseteq \ab$. Thus, $\text{Ker}\,S=\ab$ and
$\text{dim}(\ab)=36-rank(S)$.
\end{Proposition}\noin \textsc{Proof.} Straightforward computations show that an action by $e_1$ on any vector
is always a multiple of $w_{12}$ and an action by $e_2$ on any
vector is always a multiple of $\wh{w_{23}}$. It follows that an
action by $X_{st}$ with $1\leq s<t\leq 6$ on any vector is always
proportional to $w_{st}$ and an action by $\wh{X_{st}}$ with $2\leq
s<t\leq 6$ on any vector is always proportional to $\wh{w_{st}}$.
Further, by definition of the $X$'s which carry a triangle, an
action by $\ovs{\tp}{X_{st}}$ (resp $\ovs{\tpp}{X_{st}}$) with
$s\geq 3$ (resp $s\geq 4$) on any vector is always proportional to
$\ovs{\tp}{w_{st}}$ (resp $\ovs{\tpp}{w_{st}}$) ; an action by
$\ovs{\tppp}{X_{56}}$ (resp $\ovs{\triangle^l}{X_{56}}$) on any
vector always yiels a multiple of
$\ta$ (resp $\lr$). \\
In other words, each line of the matrix $S$ corresponds to the
action of exactly one conjugate from Definition $2$.\\

\noin We computed with Maple the rank of the sum matrix $S$ for the
values of $l$ and $r$ that annihilate the determinant. We found the
following values.
$$\begin{array}{ccc}
\text{When}&l=r^3&rk(S)=21\\
\text{When}&l=-r^3&rk(S)=6\\
\text{When}&l=\unsur{r^3}&rk(S)=16\\
\text{When}&l=-\unsur{r^9}&rk(S)=30\\
\text{When}&l=\unsur{r^{21}}&rk(S)=35\\
\end{array}$$

\noin The respective dimensions above correspond to the exponents in
the determinant. In particular, for each of these values of $l$ and
$r$, the module $\ab$ is non-trivial. Since it is not the whole
space $V_6$ either, this shows that for these values the
representation $\n^{(6)}$ is reducible. The necessary and sufficient
conditions are gathered in the following proposition.
\begin{Proposition}
The representation $\n^{(6)}$ is reducible if and only if
$$l\in\bigg\lb r^3,-r^3,\unsur{r^3},-\unsur{r^9},\unsur{r^{21}}\bigg\rb$$
\end{Proposition}
\noin We conclude the proof of Theorem $1$ by noticing that
$$\n^{(6)}(e_ie_j)=0\qquad\text{when $i\not\sim j$}$$
This is simply because our representation $\n^{(6)}$ is built inside
the $B(E_6)$-module $M\otimes_{\mathcal{H}}\,F$ that was introduced
earlier. Then, $\n^{(6)}$ is an irreducible representation of
$$B(E_6)e_1B(E_6)/<B(E_6)e_ie_j\,B(E_6)>_{i\not\sim j}$$
This quotient of ideals is called $I_1/I_2$ in \cite{CGW}. As part
of their work in \cite{CGW}, for each irreducible representation of
the Hecke algebra of type $A_5$ of degree $d$, Arjeh Cohen, Di\'e
Gijsbers and David Wales build an irreducible representation of
$I_1/I_2$ of degree $36\,d$. They show all the inequivalent
irreducible representations of $I_1/I_2$ are obtained this way.
There are only two inequivalent representations of the Hecke algebra
of type $A_5$ of degree $1$. Thus, there are only two inequivalent
irreducible representations of $I_1/I_2$ of degree $36$. The
representation that we built in this paper is one of them. Our $r$
is the $\unsurr$ of \cite{CGW}. Moreover, up to the change of
parameters
$$t=\unsur{lr^3}$$
and up to some rescaling of the generators, the representation of
\cite{CGW} is equivalent, as a representation of the Artin group of
type $E_6$, to the representation of \cite{CW} that was defined by
the authors and used by them to show the linearity of the Artin
group. We called the latter representation the Cohen-Wales
representation of type $E_6$. We get Theorem $1$.\\
A consequence of Lemma $2$ is that when the representation is
reducible, it is indecomposable. Then the BMW algebra is not
semisimple. So, for the values of $l$ and $r$ of Proposition $5$,
the algebra is not semisimple. Since $r$ and $-\unsurr$ play
identical roles, we obtain a second set of values of the parameters
for which the algebra is not semisimple, as in the statement of
Theorem $2$.\\\\
We conclude this paper by noting that the method that we use allows
to find Hecke algebra representations of type $E_6$.\\\\
\noin\textbf{Acknowledgments.} The author thanks David Wales for
valuable remarks and comments during the preparation of the
manuscript.

\end{document}